\documentclass{amsart}
\usepackage{bm}
\usepackage{amsmath,amssymb, amscd}
\usepackage{graphicx}

\newtheorem{theorem}{Theorem}[section]
\newtheorem{lemma}[theorem]{Lemma}
\newtheorem{proposition}[theorem]{Proposition}
\newtheorem{corollary}[theorem]{Corollary}

\theoremstyle{definition}
\newtheorem{definition}[theorem]{Definition}

\numberwithin{equation}{section}

\theoremstyle{remark}
\newtheorem{remark}[theorem]{Remark}

\begin{document}
\title[Local Gromov--Witten invariants of cubic surfaces]{Local Gromov--Witten invariants of cubic surfaces via nef toric degeneration}
\author{Yukiko Konishi}
\address{Graduate School of Mathematical Sciences, The University of Tokyo,
 3-8-1 Komaba, Meguro, Tokyo 153-8914 Japan}
\email{konishi@ms.u-tokyo.ac.jp}
\author{Satoshi Minabe}
\address{Graduate School of Mathematics, Nagoya University,
Nagoya 464-8602, Japan}
\email{minabe@yukawa.kyoto-u.ac.jp}

\begin{abstract}
We compute local Gromov--Witten invariants of cubic surfaces at all genera. 
We use a deformation of a cubic surface to a nef toric surface and  
the deformation invariance of Gromov--Witten invariants. 
\end{abstract}
\maketitle

\section{Introduction}
A del Pezzo surface $S_d$ of degree $d$ ($1\leq d\leq 9$)
\footnote{In physics literatures, $S_d$ is usually denoted by $dP_{9-d}$ or $B_{9-d}$. 
Here we follow the notation used in 
\cite[\S 0]{MukaiMori}.
A brief account of the classification of 
del Pezzo surfaces can be found there.} 
is a smooth surface%
\footnote{In this article, a surface means 
an algebraic surface over $\Bbb{C}$.} 
whose anticanonical divisor $-K_{S_d}$ is ample
and $(-K_{S_d})^2=d$.
For a smooth projective surface $X$, 
the local Gromov--Witten (GW) invariant 
is  a rational number  defined by the integral of
a certain class,
which is determined by the canonical divisor $K_X$,
on the moduli stack of stable maps to $X$ 
\cite{CKYZ, KlemmZaslow}.
Local GW invariants of del Pezzo surfaces 
have been intensively studied
in physics in relation to the  non-critical string
by various methods:
mirror symmetry, 
Seiberg--Witten curve technique and so on 
(see e.g. \cite{LMW}).
In the case of toric del Pezzo surfaces (i.e. $6\leq d\leq 9$),
a powerful method based on 
the duality to the Chern--Simons theory 
enables us to write down an explicit formula for the generating function
at all genera
\cite{DFG1, DFG2, AMV, Iqbal}.
The formula  was proved in \cite{Zhou} based on 
the virtual localization \cite{Kontsevich,GrPa}
together with a formula for Hodge integrals \cite{LLZ}.
In a recent interesting work \cite{DS}, 
Diaconescu and Florea 
proposed a closed formula for the generating function of
nontoric del Pezzo surfaces $S_i$ ($1\leq i\leq 5$)
for all genera  by 
using the conjectural ruled vertex formalism \cite{DFS}.

Our  modest goal is to obtain a formula for the generating function of
local GW invariants of $S_3$ at all genera.
$S_3$ is isomorphic to $\mathbb{P}^2$ blown-up at 6 points in a general 
position and it is also realized as a smooth cubic surface in $\mathbb{P}^3$.
It is not toric but have a (unique) smooth nef toric degeneration $S_3^0$
(a smooth toric surface with the nef anticanonical divisor 
which is deformation equivalent to $S_3$).
A main idea is to use 
the deformation invariance of local GW invariants 
as in  \cite{DS,Ueda}
and reduce the computation to those of $S_3^0$
where we can apply the virtual localization.
Here we remark that our results are limited to $S_k$ ($k=3,4,5$)
since $S_1$ and $S_2$ do not admit nef toric degenerations.

The results of this paper are as follows.
We first prove that in the case of a smooth projective 
surface with the nef anticanonical divisor,
local GW invariants
are equal to ordinary GW invariants of a projective bundle compactification
of the total space of the 
canonical line bundle (Proposition \ref{prop:local-global}).
Our proof is based on  the
virtual localization with respect to the 
$\mathbb{C}^*$-action in the fiber direction.
Then the  deformation invariance of the latter \cite{LiTian,Siebert}
implies that of the former (Proposition \ref{prop:deformation-equivalence}).
Next we introduce the toric surface  $S_3^0$ and show that it
is the nef toric degeneration of $S_3$
(Proposition \ref{prop:toric-degeneration}).
Then we derive a formula for
the generating function of local GW invariants of $S_3^0$
by the virtual localization  
(Lemma \ref{prop:generating-function-S_3^0}).
Finally we obtain a formula for the generating function of 
local GW invariants of $S_3$ via those of $S_3^0$
by the deformation invariance (Theorem  \ref{thm:formula}).

The organization of the paper is as follows.
In Section \ref{sec:localGW},
we give a definition of local GW invariants 
and show 
the deformation invariance.
In Section \ref{sec:cubic-surfaces}, 
we summarize necessary facts about cubic surfaces $S_3$.
In Section \ref{sec:deformation},
we introduce the toric surface $S_3^0$.
For completeness,
a proof of the deformation equivalence of $S_3$ and $S_3^0$
is included in Appendices \ref{sec:nef_toric} and \ref{sec:unobstructed}. 
In Section \ref{sec:comp}, 
we give formulas for the generating functions of
local GW invariants of $S_3^0$ and $S_3$.
We have computed the formula explicitly
for $\beta\in H_2(S_3,\mathbb{Z})$ 
such that $-K_{S_3}.\beta\leq 6$.
The results are listed in Section \ref{sec:GV} and Appendix \ref{sec:F_d}.

\section*{acknowledgments}
Y.K. thanks Shinobu Hosono for valuable discussions, comments 
and helpful advice on the numerical computation.
She also thanks Florin Ambro for helping her with  the proof of 
Lemma \ref{lem:toric}. 
S.M. is grateful to Hiroaki Kanno 
for many advices and encouragements. 
He is also grateful to Yuichi Nohara for discussions 
on the proof of Proposition \ref{prop:toric-degeneration}.   
Both authors would like to  thank Bogdan Florea and Kazushi Ueda  
for useful comments on the draft. 
The research of Y.K. is supported by
JSPS Research Fellowships for Young Scientists.

The authors thanks 
Hiroshi Iritani, Hiraku Nakajima, Kaoru Ono  and K\=ota Yoshioka
for pointing out an error in
our original definition of the local GW invariants.

\section{Deformation invariance of local GW invariants}
\label{sec:localGW}
In this article,
we call a smooth projective 
surface $X$ whose anticanonical divisor $-K_X$ is nef (i.e. $-K_X.[C]\geq 0$ for
all curves $C\subset X$) a nef surface.

\newcommand{\px}{\mathbb{P}(K_X\oplus \mathcal{O}_X)}

Let $X$  be a nef surface
and
$K_X$ its canonical divisor.
For $\beta\in H_2(X,\mathbb{Z})$ and $g\in\mathbb{Z}_{\geq 0}$,
let $\bar{M}_{g,0}(X,\beta)$ (resp. $\bar{M}_{g,1}(X,\beta)$)
be the moduli stack of stable maps to $X$
of genus $g$ without marked point (resp. with one marked point)
and with the second homology class $\beta$.
Let $\pi:\bar{M}_{g,1}(X,\beta)\to \bar{M}_{g,0}(X,\beta)$ be
the forgetful map of  the marked point and 
$\mu:\bar{M}_{g,1}(X,\beta)\to X$ be the evaluation at the marked point.

\newcommand{\kx}{\displaystyle{\int_{\beta}}c_1(K_X)}

\begin{definition}\label{def:localGW}
For $g\in\mathbb{Z}_{\geq 0}$ and 
$\beta\in H_2(X,\mathbb{Z})$ such that $\kx<0$,
the local Gromov--Witten invariant $N_{g,\beta}(K_X)$ of $X$
with genus $g$ and the second homology class $\beta$ is 
\begin{equation}\notag
N_{g,\beta}(K_X)=\int_{[\bar{M}_{g,0}(X,\beta)]^{vir}}
c_{top}(R^1\pi_*\mu^*  K_X),
\end{equation}
where $c_{top}$ denotes the top Chern class which is of degree
$(1-g)(\dim X -3)-\kx$.
(This is equal to the virtual dimension
of $\bar{M}_{g,0}(X,\beta)$.)%
\footnote{
The condition $\kx<0$ and the nef condition on $X$ imply
$H^0(C,f^*K_X)=0$ for $(f,C)\in \bar{M}_{g,0}(X,\beta)$.
}
\end{definition}

Let $\px$ be the projectivization of the total space of 
the vector bundle $K_X\oplus \mathcal{O}_X$
(here the canonical divisor $K_X$ and the structure sheaf $\mathcal{O}_X$
are regarded as line bundles).
This is a $\mathbb{P}^1$-bundle over $X$.
Let $\iota:X\hookrightarrow \px$ be the inclusion as the zero section of 
$K_X\subset \px$.
We define the (ordinary) GW invariant
$N_{g,\iota_*\beta}(\px)$ of 
$\px$ of genus $g$ and the second homology class $\iota_*\beta$ by 
\begin{equation}\notag
N_{g,\iota_*\beta}(\px)=\int_{[\bar{M}_{g,0}(\px,\iota_*\beta)]^{vir}} 1~.
\end{equation}
We note that
the deformation invariance is established 
for this ordinary GW invariant
\cite{LiTian,Siebert}.

\begin{proposition}\label{prop:local-global}
Let $X$ be a nef surface, $\iota:X\hookrightarrow \px$
be the inclusion as the zero section of $K_X$.
For $g\in\mathbb{Z}_{\geq 0}$ and 
$\beta\in H_2(X,\mathbb{Z})$ such that $\kx<0$,
\begin{equation}\notag
N_{g,\beta}(K_X)=N_{g,\iota_*\beta}(\px).
\end{equation}
\end{proposition}

Consider the natural
$\mathbb{C}^*$ action on $\px$ 
as the scalar multiplication in the $\mathbb{P}^1$-fiber direction.
The action induces an action on 
$\bar{M}_{g,0}(\px,\iota_*\beta)$
by moving the image curves of stable maps.
First we  show the following lemma.
\begin{lemma}\label{prop:fixed-loci}
Let $X$ be a nef surface, $\iota:X\hookrightarrow \px$
be the inclusion as the zero section of $K_X$. 
Let $\beta\in H_2(X,\mathbb{Z})$ be a class satisfying $\kx<0$.
If a stable map $(f,C)\in \bar{M}_{g,0}(\px,\iota_*\beta)$,
where $C$ is a connected curve of genus $g$ and 
$f:C\to X$ a morphism such that $[f(C)]=\iota_*\beta$, 
is fixed by the $\mathbb{C}^*$-action,
then the image $f(C)$ is contained in the zero section $\iota(X)$.
\end{lemma}

\begin{proof}
Denote the $\mathbb{P}^1$-fibration $\px\to X$ by $p$,
and let $P=[p^{-1}(a)]\in H_2(\px,\mathbb{Z})$ be the class of 
the fiber $\mathbb{P}^1$ where $a\in X$ is any point.
Let $\iota^{\infty}:X\hookrightarrow \px$ be the inclusion  as the 
zero section of $\mathcal{O}_X$ 
(the section at the infinity 
of the $\Bbb{P}^1$-bundle compactification of $K_X$). 
Note that for any 
$\alpha\in H_2(X,\mathbb{Z})$, we have
\begin{equation}\label{eq:infty}
\iota_*^{\infty}\alpha=\iota_*\alpha- \left( \int_{\alpha}c_1(K_X) \right) P ~~.
\end{equation}

Let $\gamma\in H_2(\px,\mathbb{Z})$.
If a stable map 
$(f,C)\in\bar{M}_{g,0}(\px,\gamma) $ is fixed  by the $\mathbb{C}^*$-action,
then the image of an irreducible component $C_i$ of $C$ must be
either one of these:
(i) $f(C_i)\subset \iota(X)$,
(ii) $f(C_i) \subset \iota^{\infty}(X)$ or
(iii) $f(C_i)=p^{-1}(a_i)$ $(a_i\in X)$ and $C_i\cong \mathbb{P}^1$.
So assume  that 
irreducible components $C_1,\ldots,C_k$ of $C$  are of type (i) with
$[f(C_i)]=\beta_i\in H_2(X,\mathbb{Z})$,
$C_{k+1},\ldots,C_{r}$ are of type (ii) with
$[f(C_i)]=\beta_i\in H_2(X,\mathbb{Z})$,
and that
$C_{r+1},\ldots,C_s$ are of type (iii)
with $f:C_i\to p^{-1}(a_i)$ the $d_i$-fold coverings.
Then $[f(C)]=\gamma$ is equivalent to
\begin{equation}\notag
\begin{split}
\gamma&=
\sum_{i=1}^k \iota_*\beta_i+\sum_{i=k+1}^r \iota^{\infty}_* \beta_i
+\sum_{i=r+1}^s d_i P=\sum_{i=1}^r \iota_*\beta_i
+\Big(\sum_{i=r+1}^s d_i -\sum_{i=k+1}^r \int_{\beta_i}c_1(K_X)\Big)P ~~.
\end{split}
\end{equation}
Now take $\gamma=\iota_*\beta$ with $\beta\in H_2(X,\mathbb{Z})$
satisfying $\kx<0$
and solve the above equation.
The assumption that $X$ is nef implies 
that the coefficient of $P$ in the last line is always nonnegative.
Therefore it is zero 
if and only if there is no irreducible components of type (iii)
and $\displaystyle{\int_{\beta_i}}c_1(K_X)=0$ for those of type (ii).
Then connectedness of the domain curve $C$ implies 
either $f(C)\subset \iota(X)$ or $f(C)\subset \iota^{\infty}(X)$.
For the latter case, $\displaystyle{\int_{[f(C)]}}c_1(K_X)=0$ and 
this contradicts the assumption $\kx<0$. 
Thus $f(C)\subset \iota(X)$. 
\end{proof}

\begin{proof} {\it (of Proposition \ref{prop:local-global}.)}
By Lemma \ref{prop:fixed-loci}, the $\mathbb{C}^*$-fixed 
point set is isomorphic to
$\bar{M}_{g,0}(X,\beta)$.
Then, by the virtual localization \cite{GrPa},
\begin{equation}\notag
N_{g,\iota_*\beta}(\px)=
\int_{[\bar{M}_{g,0}(X,\beta)]^{vir}}
e_{\mathbb{C}^*}(R^{1}\pi_*\mu^*  K_X).
\end{equation}
Here $e_{\mathbb{C}^*}$ is the equivariant Euler class.
(In the equation below \cite[(24)]{GrPa},
the nontrivial contribution comes only from  the factor  $e(B_5^m)$;
$e(B_2^m)$ does not contribute because $\kx<0$.)
Since the LHS is independent of the weight, so is the RHS and
we can replace it  with the nonequivariant integral.
\end{proof}

\begin{proposition}\label{prop:deformation-equivalence}
Let $X$ be a nef surface and
$X'$ be a nef surface which is deformation equivalent to $X$.
Let
$\beta\in H_2(X,\mathbb{Z})$ be a class satisfying  $\kx<0$
and $\beta'\in H_2(X',\mathbb{Z})$ be the class corresponding to $\beta$ 
under a deformation.
Then
$N_{g,\beta}(K_X)=N_{g,\beta'}(K_{X'})$
for $g\in \mathbb{Z}_{\geq 0}$.
\end{proposition}

\begin{proof}
Since $X$ and $X'$ are deformation equivalent,
$\px$ and $\mathbb{P}(K_{X'}\oplus \mathcal{O}_{X'})$ are also deformation
equivalent.
Let $\iota:X\hookrightarrow \px$ and $\iota':X'\hookrightarrow 
\mathbb{P}(K_{X'}\oplus \mathcal{O}_{X'})$
be the inclusions as the zero sections of $K_X$ and $K_{X'}$
respectively.

We have
\begin{equation}\notag
N_{g,\beta}(K_X)=N_{g,\iota_*\beta}(\px)=
N_{g,\iota_*'\beta'}(\mathbb{P}(K_{X'}\oplus \mathcal{O}_{X'}) )
=N_{g,\beta'}(K_{X'}).
\end{equation}
The middle equality follows from
the deformation invariance of ordinary GW invariants \cite{LiTian,Siebert}.
The first and the third equalities follow from Proposition
\ref{prop:local-global}.
\end{proof}

\section{Cubic surfaces $S_3$}
\label{sec:cubic-surfaces}
Here we summarize some facts on cubic surfaces, 
see e.g.  \cite[Ch. V, 4]{Ha} for details. 

Let $S_3$ be a cubic surface.  
$S_3$ is realized as a blowing up $\pi: S_3 \to \Bbb{P}^2$ 
at six points in a general position.  
Let $e_1, \cdots, e_6$ be the classes of 
the exceptional curves of $\pi$ and 
$l$ be 
the class of a line in $\Bbb{P}^2$ pulled back by $\pi$. 
Then $l, e_1, \cdots, e_6$ is a basis of $\textrm{Pic}(S_3)$.
Their intersections are
\begin{equation*}
l^2=1,~~~ e_i^2=-1,~~~ l.e_i = 0,~~~ e_i.e_j = 0 ~\mathrm{if}~ i \neq j. 
\end{equation*} 
Let $h$ be the class of the hyperplane section of $\Bbb{P}^3$. 
Then we have  
\begin{equation*}
h =  -K_{S_3} = 3l - \sum_{i=1}^6 e_i~.
\end{equation*} 

It is a classical fact that $S_3$ contains exactly twenty-seven 
lines which are given as follows:
\begin{equation*}
e_i~(i=1,\cdots, 6), \quad  l-e_i-e_j~(1\leq i < j\leq 6), 
 \quad 2l-\sum_{i\neq j}e_i~~(j=1, \cdots, 6).
\end{equation*}
Each one of these is a exceptional curve of the first kind.
These twenty-seven lines are 
the minimal generators of 
the Mori cone (the cone 
generated by effective divisors on 
$X$ modulo numerical equivalence) (cf. \cite[(0.6)]{MukaiMori}).

It is well-known that the Weyl group $W_{E_6}$ of type $E_6$ acts on 
$\textrm{Pic}(S_3)$ 
as symmetries of configurations of twenty seven lines. 
Its generators are given as follows.  
\begin{equation}\notag
\begin{split}
&s_i:e_i\leftrightarrow e_{i+1} ~~(1\leq i\leq 5),
\\
&s_6:e_1\mapsto l-e_2-e_3,\quad
     e_2\mapsto l-e_1-e_3,\quad
     e_3\mapsto l-e_1-e_2,\quad
     l\mapsto 2l-e_1-e_2-e_3.
\end{split}
\end{equation}
It is known  \cite[\S 4]{Dolgachev} that $W_{E_6}$ coincides with  
the group of automorphisms of $\text{Pic}(S_3)$ 
which preserve the intersection form, the canonical class, 
and the semigroup of effective classes. 
It is also known that such an automorphism on $\text{Pic}(S_3)$ 
comes from an isomorphism of $S_3$.  

Hereafter we identify $\textrm{Pic}(S_3)$ with 
$H^2(S_3,\mathbb{Z})\cong H_2(S_3,\mathbb{Z})$.
The next lemma was shown in \cite[\S 2.4]{Hosono}. 
\begin{lemma}\label{prop:weyl-invariance}
$N_{g,\beta}(K_{S_3})=N_{g,w(\beta)}(K_{S_3})$
for  $w\in W_{E_6}$.
\end{lemma}
\begin{proof}
Since the action of $w$ on $H_2(S_3, \Bbb{Z})$ 
is induced from an isomorphism $\psi : S_3 \to S_3$,  
we have $N_{g ,\beta} (K_{S_3}) = 
N_{g, \psi_{*}\beta}(K_{S_3}) = N_{g, w(\beta)}(K_{S_3})$.
\end{proof}
\section{Nef toric surfaces deformation equivalent to 
 $S_3$, $S_4$, and $S_5$}
\label{sec:deformation}

\begin{figure}[t]
\unitlength .10cm
\begin{center}
\begin{picture}(30,30)(40,-10)
\thicklines
\put(0,0){\line(1,0){10}}
\put(0,0){\line(0,1){10}}
\put(0,0){\line(-1,2){10}}
\put(0,0){\line(-1,1){10}}
\put(0,0){\line(-1,0){10}}
\put(0,0){\line(-1,-1){10}}
\put(0,0){\line(0,-1){10}}
\put(0,0){\line(1,-1){10}}
\put(0,0){\line(2, -1){20}}
\put(20,-10){\line(-1,1){30}}
\put(20,-10){\line(-1,0){30}}
\put(-10,20){\line(0,-1){30}}
\put(11,0){$v_1$}
\put(0,10){$v_2$}
\put(-14,19){$v_3$}
\put(-14,9){$v_4$}
\put(-14,-1){$v_5$}
\put(-14,-10){$v_6$}
\put(0,-12){$v_7$}
\put(10,-12){$v_8$}
\put(20,-12){$v_9$}
\put(14,0){\vector(1,0){10}}
\put(40,0){\line(1,0){10}}
\put(40,0){\line(0,1){10}}
\put(40,0){\line(-1,1){10}}
\put(40,0){\line(-1,0){10}}
\put(40,0){\line(-1,-1){10}}
\put(40,0){\line(0,-1){10}}
\put(40,0){\line(1,-1){10}}
\put(40,0){\line(2, -1){20}}
\put(60,-10){\line(-1,1){20}}
\put(60,-10){\line(-1,0){30}}
\put(30,10){\line(0,-1){20}}
\put(30,10){\line(1,0){10}}
\put(54,0){\vector(1,0){10}}
\put(80,0){\line(1,0){10}}
\put(80,0){\line(0,1){10}}
\put(80,0){\line(-1,1){10}}
\put(80,0){\line(-1,0){10}}
\put(80,0){\line(0,-1){10}}
\put(80,0){\line(1,-1){10}}
\put(90,0){\line(0,-1){10}}
\put(90,0){\line(-1,1){10}}
\put(90,-10){\line(-1,0){20}}
\put(70,10){\line(0,-1){20}}
\put(70,10){\line(1,0){10}}
\put(80,0){\line(-1,-1){10}}
\put(94,0){\vector(1,0){10}}
\put(120, 0){\line(1,0){10}}
\put(120,0){\line(0,1){10}}
\put(126,-4){$p_2$}
\put(120,0){\line(-1,1){10}}
\put(116, 6){$p_3$}
\put(120,0){\line(-1,0){10}}
\put(116, -4){$p_1$}
\put(120,0){\line(0,-1){10}}
\put(120,0){\line(1,-1){10}}
\put(120,10){\line(1,-1){10}}
\put(120,10){\line(-1,0){10}}
\put(110,0){\line(0,1){10}}
\put(110,0){\line(1,-1){10}}
\put(130, -10){\line(0,1){10}}
\put(130, -10){\line(-1,0){10}}
\end{picture}
\end{center}
\caption{$S_3^0 \to S_4^0 \to S_5^0 \to S_6$.
}\label{fig:fanS}
\end{figure}
Let $S_3^0$, $S_4^0$, and $S_5^0$ be the nef toric surfaces whose 
fans are given in Figure \ref{fig:fanS}. 
Here the nine one-dimensional cones  of $S_3^0$ are generated by
\begin{equation}\notag 
\begin{split}
&v_1=(1,0),\quad v_2=(0,1),\quad v_3=(-1,2),\quad
v_4=(-1,1),\quad v_5=(-1,0),\\
&\quad v_6=(-1,-1),\quad
v_7=(0,-1),\quad v_8=(1,-1),\quad v_9=(2,-1).
\end{split}
\end{equation}
Let 
the fan of the toric del Pezzo surface $S_6$ be given 
in Figure \ref{fig:fanS}
and
let $p_1,p_2,p_3$ be the torus fixed points of $S_6$ corresponding to
the two-dimensional cones generated by $(v_5,v_7),(v_8,v_1),(v_2,v_4)$.
$S_3^0$ (resp. $S_4^0,S_5^0$) is obtained 
by blowing up $S_6$ at $p_1,p_2,p_3$ (resp. $p_1,p_2$ and $p_1$).
$S_k^0$ contains $(-2)$-curves and its anticanonincal divisor is nef but
not ample.   

\begin{proposition}\label{prop:toric-degeneration}
$S_k^0$ ($k=3,4,5$) is deformation equivalent to $S_k$. 
\end{proposition}
A proof will be given 
in Appendix \ref{sec:nef_toric} (see Proposition \ref{prop:def_type}).

Now let us explain the geometry of the nef toric surface $S_3^0$.
The torus-invariant divisors  $C_i$ ($1\leq i\leq 9$) corresponding to
$v_i$ have the intersections:
\begin{equation}\label{eq:intersection-S_3^0}
C_i.C_{i+1}=1,\quad
C_i.C_j=0\,\,(j\neq i,i\pm 1),\quad
C_i^2=\Bigg\{\begin{array}{cc}-1 &(i=3,6,9)~,\\
                             -2 &(i=1,2,4,5,7,8)~,\end{array}
\end{equation}
and the canonical divisor $K_{S_3^0}$ is rationally equivalent to 
$-C_1-\cdots -C_9$.
The Mori cone is generated by
$C_1,\ldots,C_9$ \cite[Proposition 2.26]{Oda}.

Note that $\textrm{Pic}(S_3^0) \cong \textrm{Pic}(S_3)$ 
and an isomorphism is given by the following. 
\begin{equation}\label{eq:identification}
\begin{split}
C_1\mapsto e_2-e_5, \qquad
C_2\mapsto l-e_2-e_3-e_6,\qquad   
C_3\mapsto  e_6,\\
C_4\mapsto e_3-e_6,\qquad
C_5\mapsto l-e_1-e_3-e_4,\qquad
C_6\mapsto e_4,\\
C_7\mapsto  e_1-e_4,\qquad 
C_8\mapsto l-e_1-e_2-e_5,\qquad   
C_9\mapsto  e_5.
\end{split}
\end{equation}
This 
is explained as follows.
First, in $S_6$, we regard
the torus-invariant divisors $C_1',C_4',C_7'$
corresponding to $v_1,v_4,v_7$
as the exceptional curves of blowing up of $\mathbb{P}^2$
and identify them  with $e_2,e_3,e_1$.
The torus-invariant divisors $C_2',C_5',C_8'$ corresponding to $v_2,v_5,v_8$
are identified with 
the proper transforms $l-e_2-e_3,l-e_1-e_3,l-e_1-e_2$ 
of lines in $\mathbb{P}^2$. Then in $S_3^0$, 
$C_3,C_6,C_9$ are exceptional curves of the blowup at
$p_3,p_1,p_2$ and we identify them with $e_6,e_4,e_5$.
For $i=1,2,4,5,7,8$,
$C_i$ is the proper transform of $C_i'$.
(This identification 
can be seen from the construction of
a deformation in the proof of Proposition \ref{prop:def_type}.)

From here on, we identify $\rm{Pic}(S_3^0)$ with 
$H^2(S_3^0,\mathbb{Z})\cong H_2(S_3^0,\mathbb{Z})$.
\begin{theorem}\label{thm:1}
For $g\in \mathbb{Z}_{\geq 0}$ and
$\beta\in H_2(S_3,\mathbb{Z})$ such that $K_{S_3}.\beta<0$,
$$N_{g,\beta}(K_{S_3})=N_{g,\beta'}(K_{S_3^0})~,$$ 
where
$\beta'\in H_2(S_3^0,\mathbb{Z})$ is 
the class corresponding to $\beta$ by eq. (\ref{eq:identification}).
\end{theorem}

\begin{proof}
This follows from
Propositions 
\ref{prop:deformation-equivalence} and
\ref{prop:toric-degeneration}.
\end{proof}

\begin{remark}\label{rem:4and5}
The statements similar to Theorem \ref{thm:1} hold for $S_4,S_5$:
local GW invariants of
$S_4$ and $S_5$ are the same as those of $S_4^0$ and $S_5^0$.  
Their generating functions also have expressions 
analogous to the formula 
for $S_3$ (which will be stated in Theorem \ref{thm:formula}).
Local GW invariants of $S_4$ and $S_5$ appear among those of $S_3$
with a natural identification
of second homology classes
$H_2(S_3,\mathbb{Z})=H_2(S_4,\mathbb{Z})\oplus \mathbb{Z}e_6=
H_2(S_5,\mathbb{Z})\oplus \mathbb{Z}e_5\oplus \mathbb{Z}e_6$.
See \cite[\S 6]{KonishiMinabe}. 
\end{remark}

\section{Formula for the generating function of local GW invariants of $S_3$}
\label{sec:comp}
\subsection{}
First we consider the generating function of local GW invariants 
of $S_3^0$ with $\beta\in H_2(S_3^0,\mathbb{Z})$ such that
$K_{S_3^0}.\beta<0$.
Take a basis $c_1,\ldots,c_7$ of $H_2(S_3^0,\mathbb{Z})$
and let $X_1,\ldots,X_7$ be associated formal variables.
For $\beta=a_1c_1+\dots+a_7 c_7 \in H_2(S_3^0,\mathbb{Z})$, 
denote $X_1^{a_1}\dots X_7^{a_7}$ by $X^{\beta}$.
We write the generating function as
\begin{equation}\notag
F_{S_3^0}=\sum_{\begin{subarray}{c}
\beta\in H_2(S_3^0,\mathbb{Z}),\\K_{S_3^0}.\beta<0
\end{subarray}}
\sum_{g\geq 0} N_{g,\beta}(K_{S_3^0})
\lambda^{2g-2}X^{\beta}.
\end{equation}
Let $t_i=X^{[C_i]}$ ($1\leq i\leq 9$) and
$s_i=C_i^2$ (See (\ref{eq:intersection-S_3^0})).
Define $Z_{S_3^0}$ by
\begin{equation}\notag
Z_{S_3^0}=\prod_{i=1}^9 \sum_{\nu^i} 
((-1)^{s_i}t_i)^{|\nu^i|}
e^{\sqrt{-1}\lambda s_i\frac{\kappa(\nu^i)}{2}}
W_{\nu^i,\nu^{i+1}}(e^{\sqrt{-1}\lambda}).
\end{equation}
Here each $\nu^i$ ($1\leq i \leq 9$) runs over the set of partitions
and $\nu^{10}=\nu^1$ is assumed.
For partitions 
$\mu=(\mu_1,\mu_2,\ldots)$  and $\nu=(\nu_1,\nu_2,\ldots)$,
\begin{equation}\notag
W_{\mu,\nu}(q)=s_{\mu}(q^{\rho})s_{\nu}(q^{\mu+\rho}) 
\in \mathbb{Q}(q^{\frac{1}{2}}),\quad
|\mu|=\sum_{i\geq 1}\mu_i, \quad
\kappa(\mu)=\sum_{i\geq 1}\mu_i(\mu_i-2i+1),
\end{equation}
where 
$q^{\mu+\rho}=(q^{\mu_i-i+\frac{1}{2}})_{i\geq 1}$, 
$q^{\rho}=(q^{-i+\frac{1}{2}})_{i\geq 1}$ and
$s_{\mu}$ denotes the Schur function.
Define $Z_{(-2)}(t)$ by
\begin{equation}\notag
Z_{(-2)}(t)=
\exp\Bigg[
-\sum_{j\geq 1}\frac{1}{j}\Big(2\sin\frac{j\lambda}{2}\Big)^{-2}t^j 
\Bigg]~.
\end{equation}

\begin{lemma}\label{prop:generating-function-S_3^0}
\begin{equation}\notag
\exp\big(F_{S_3^0}\big)=
\frac{Z_{S_3^0}}
{\prod_{i=1,4,7}Z_{(-2)}(t_i)Z_{(-2)}(t_{i+1})Z_{(-2)}(t_it_{i+1})}~~.
\end{equation}
\end{lemma}

\begin{proof}
Recall that $S_3^0$ has a canonical
$T=(\mathbb{C}^*)^2$-action determined by its fan.
Let $K_{S_3^0}^T=-C_1-\ldots-C_9$ be an $T$-invariant divisor.
For any $\beta\in H_2(S_3^0,\mathbb{Z})$ 
and $g\in \mathbb{Z}_{\geq 0}$, 
define $N_{g,\beta}^T(S_3^0)$ by the following equivariant integral:
\begin{equation} \notag 
N_{g,\beta}^T(S_3^0)=
\int_{[\bar{M}_{g,0}(S_3^0,\beta)^T]^{vir} }
\frac{e_T(R^{1}\pi_*\mu^* K_{S_3^0}^T)}
{e_T(R^0\pi_*\mu^* K_{S_3^0}^T)}
\frac{1}{e_T(Norm)}~~.
\end{equation}
Here
$\bar{M}_{g,0}(S_3^0,\beta)^T$ is the fixed point set
of the induced $T$-action,
$e_T$ denotes the equivariant Euler class and
$Norm$ is the virtual normal bundle determined by the
obstruction theory \cite[eqs. (23)(24)]{GrPa}.
Note that 
$N_{g,\beta}^T(S_3^0)=0$
if there is no effective divisors of the form
 $\sum_{1\leq i\leq 9} a_i[C_i]$ ($a_i\in\mathbb{Z}_{\geq 0}$)
which are rationally equivalent to $\beta$
because $\bar{M}_{g,0}(S_3^0,\beta)^T$ is empty.

Consider the exponential of the 
generating function for {\em all} classes
\begin{equation}\label{eq:generating-function2}
\exp\Bigg[ \sum_{\beta\in H_2(S_3^0,\mathbb{Z})}
\sum_{g\geq 0} N_{g,\beta}^T(S_3^0)\lambda^{2g-2}X^{\beta}
\Bigg]~~.
\end{equation}
Carrying  out the localization calculation
in the same way as \cite{Zhou}%
\footnote{The contribution to $N_{g,\beta}^T(S_3^0)$ from a fixed locus
turns out to be completely the same as \cite[eqs. (13)(16)]{Zhou}.
Thus the summation over genera, second homology classes and fixed loci
proceeds in the same manner.}
and
using the formula for Hodge integrals \cite[Theorem 1]{LLZ},
we see that $(\ref{eq:generating-function2})$ is equal to $Z_{S_3^0}$.


Next we have to subtract the contributions coming from 
classes $\beta$ which does not satisfy
$K_{S_3^0}.\beta<0$.
Note that such effective classes
are of the forms $a[C_1]+b[C_2]$, $a[C_4]+b[C_5]$ or 
$a[C_7]+b[C_8]$
($a,b\in\mathbb{Z}_{\geq 0}$).
Therefore
\begin{equation}\label{eq:-2contribution}
\exp
\Bigg[
\sum_{\begin{subarray}{c}\beta\in H_2(S_3^0,\mathbb{Z})\\
           K_{S_3^0}.\beta\geq 0\end{subarray}}
\sum_{g\geq 0} N_{g,\beta}^T(S_3^0)\lambda^{2g-2}X^{\beta}
\Bigg]
=
\prod_{i=1,4,7}\exp\Bigg[
\sum_{a,b\in \mathbb{Z}_{\geq 0}}
\sum_{g\geq 0} N_{g,a[C_i]+b[C_{i+1}]}^T(S_3^0)\lambda^{2g-2}t_i^at_{i+1}^b
\Bigg].
\end{equation}
The $i=1$ factor is easily obtained by
setting $t_3=t_4=\dots=t_9=0$ in (\ref{eq:generating-function2}).
It is equal to 
$$
Z_{S_3^0}|_{t_3=t_4=\dots=t_9=0}=Z_{(-2)}(t_1)Z_{(-2)}(t_2)Z_{(-2)}(t_1t_2)~~
.$$
The $i=4,7$ factors are similar.
Dividing (\ref{eq:generating-function2}) by (\ref{eq:-2contribution}),
we obtain
\begin{equation}\notag
\exp\Bigg[ 
\sum_{\begin{subarray}{c}\beta\in H_2(S_3^0,\mathbb{Z}),\\
K_{S_3^0}.\beta<0\end{subarray}}
\sum_{g\geq 0} N_{g,\beta}^T(S_3^0)\lambda^{2g-2}X^{\beta}\Bigg]
=\frac{Z_{S_3^0}}
{\prod_{i=1,4,7}Z_{(-2)}(t_i)Z_{(-2)}(t_{i+1})Z_{(-2)}(t_it_{i+1}) }.
\end{equation}
By the 
virtual localization \cite{GrPa},
$N_{g,\beta}^T(S_3^0)=N_{g,\beta}(S_3^0)$ 
for $\beta$ such that $K_{S_3^0}.\beta<0$.
Thus We complete our proof.
\end{proof}
\subsection{}
Next we study the generating function of local GW invariants of  $S_3$. 
Let $Q=(Q_1,\ldots,Q_6,Q_7)$ be a set of formal variables 
and denote $Q_1^{a_1}Q_2^{a_2}\ldots Q_7^{a_7}$ by $Q^{\beta}$
for $\beta=a_1e_1+\cdots+a_6 e_6+a_7 l\in H_2(S_3,\mathbb{Z})$.
Define
\begin{equation}\notag
\begin{split}
&F_d= 
\sum_{\begin{subarray}{c}\beta\in H_2(S_3,\mathbb{Z}),\\
                       -K_{S_3}.\beta=d\end{subarray}}
\sum_{g\in\mathbb{Z}_{\geq 0}}
N_{g,\beta}(K_{S_3})\lambda^{2g-2} Q^{\beta}~, 
\qquad (d\in\mathbb{Z}_{\geq 1}),
\end{split}
\end{equation}
and $F_{S_3}: =\displaystyle{\sum_{d\geq 1}} F_d$.

\begin{theorem} \label{thm:formula}
With the following identification of the parameters
\begin{equation}\label{eq:id-parameters}
\begin{split}
&t_1=Q^{e_2-e_5},\quad
t_2=Q^{l-e_2-e_3-e_6},\quad
t_3=Q^{e_6},\quad
t_4=Q^{e_3-e_6},\quad 
t_5=Q^{l-e_1-e_3-e_4},
\\
&t_6=Q^{e_4},\quad
t_7=Q^{e_1-e_4},\quad
t_8=Q^{l-e_1-e_2-e_5},\quad
t_9=Q^{e_5},
\end{split}
\end{equation}
we have 
\begin{equation}\notag
\exp\big(F_{S_3}\big)=\exp\big(F_{S_3^0}\big).
\end{equation}
\end{theorem}

\begin{proof}
This follows from Theorem \ref{thm:1} 
and Lemma \ref{prop:generating-function-S_3^0}.
The identification (\ref{eq:id-parameters})
is determined by (\ref{eq:identification}).
\end{proof}

\begin{remark} 
In \cite{DS}, Diaconescu and Florea obtained a formula for $F_{S_3}$
which is different from ours
(eq. (3.14) for $k=5$ in {\it loc. cit.}).
It would be an interesting problem to show that
these two formulas are  equivalent.
\end{remark}

Define $m(\beta)$ for 
$\beta \in H_2(S_3,\mathbb{Z})$ by
\begin{equation}\notag
m(\beta)=\frac{1}{\#\{w\in W_{E_6} \mid w(\beta)=\beta \}} \sum_{w\in W_{E_6}}
Q^{w(\beta)}~.
\end{equation}
By Lemma \ref{prop:weyl-invariance}, $F_d$ should be 
written in terms of these.
$F_d$ up to $d=6$ are shown in Appendix \ref{sec:F_d}.

\section{Gopakumar--Vafa invariants}\label{sec:GV}
Let
$n_{\beta}^g(K_{S_3})$ ($g\in\mathbb{Z}_{\geq 0}$,
$\beta\in H_2(S_3,\mathbb{Z})$) be numbers defined
by the following :
\begin{equation}\notag
F_{S_3}=\sum_{\beta\in H_2(S_3,\mathbb{Z})}
\sum_{g\in\mathbb{Z}_{\geq 0}}
\sum_{k\geq 1} \frac{n_{\beta}^{g}(K_{S_3})}{k}
\Big(2\sin \frac{k \lambda}{2}\Big)^{2g-2}Q^{k\beta}.
\end{equation}
$n_{\beta}^g(K_{S_3})$ are called Gopakumar--Vafa invariants \cite{GV} .
They are listed in Table \ref{tab:GV}.

\begin{table}[ht]
\begin{equation}\notag
\begin{array}{|c|r|r|r|crrrrrr|}\hline
d&\beta &\#\mathcal{O}(\beta)&\text{genus}&g
   & 0  & 1 & 2 &  3& 4 &5\\\hline 
1& e_6                        &27  &0 && 1 &&&&&\\\hline
2&-e_1+l                      &27  &0 && -2 &&&&&\\\hline
3&l                           &72  &0 && 3  &   &&&&\\
 &-e_1-e_2-e_3-e_4-e_5-e_6+3l &1   &1 && 27 &-4 &&&&\\\hline
4&-e_1-e_2+2l                 &216 &0 &&-4  &   &&&&\\\
 &-e_1-e_2-e_3-e_4-e_5+3l     &27  &1 &&-32 & 5 &&&&\\\hline
5&-e_1+2l                     &432 &0 && 5  &   &&&&\\
 &-e_1-e_2-e_3-e_4+3l         &216 &1 && 35 &-6  &&&&\\
 &-2e_2-e_2-e_3-e_4-e_5-e_6+4l &27  &2 && 205& -68& 7 &&&\\\hline  
6&-2e_1-e_2+3l                 &432 &0 && -6  &      &     &  &&\\
 &2l                           &72  &0 && -6  &      &     &  &&\\
 &-e_1-e_2-e_3+3l              &720 &1 && -36 & 7    &     &  &&\\
 &-2e_1-e_2-e_3-e_4-e_5+4l     &270 &2 &&-198 & 72   & -8  &  &&\\
 &-e_1-e_2-e_3-e_4-e_5-e_6+4l  &72  &3 && -936& 498  &-108 & 9&&\\
 &-2e_1-2e_2-2e_3-2e_4-2e_5-2e_6+6l &1&4&&-3780& 2636&-846 &141&-10&\\\hline  
\end{array}
\end{equation}
\caption{Gopakumar--Vafa invariants $n_{\beta}^g(K_{S_3})$}

\label{tab:GV}
\end{table}

\begin{remark}
\begin{enumerate}
\item[(a)]
Gopakumar--Vafa invariants $n_{\beta}^g(K_{S_3})$ of $S_3$ are integers. 
Moreover, for each $\beta$, 
$n_{\beta}^g(K_{S_3})$ is equal to zero for all but finite $g$.  
This follows from the same statement for the toric surface 
$S_3^0$ (\cite{Peng, konishi}).

\item[(b)]
One could observe that
$n_{\beta}^g(K_{S_3})$ in Table \ref{tab:GV} are zero
if
$g$ is larger than
the genus 
$\beta.(\beta+K_{S_3})/2+1$
of a nonsingular curve which belongs to $\beta$.

\item[(c)] 
The results are in agreement with previous results in
\cite[Table 3]{MNW},
\cite[Table 1, $n=6$]{LMW},
\cite[Table 7, $X_3(1,1,1,1)$]{CKYZ}
obtained by the B-model calculation of mirror symmetry.
Also compare with \cite[Table 7]{KKV}.

\end{enumerate}
\end{remark}
\appendix
\section{Nef toric surfaces and their deformations}\label{sec:nef_toric}

\begin{figure}[t]
\unitlength .10cm
\begin{center}
\begin{picture}(30,30)(40,13)
\thicklines
\put(-18,-2){\framebox(4, 4){16}}
\put(-10,10){\line(1,0){5}}
\put(-10,10){\line(0,1){5}}
\put(-10,10){\line(-1,2){5}}
\put(-10,10){\line(-1,1){5}}
\put(-10,10){\line(-1,0){5}}
\put(-10,10){\line(-1,-1){5}}
\put(-10,10){\line(0,-1){5}}
\put(-10,10){\line(1,-1){5}}
\put(-10,10){\line(2, -1){10}}
\put(0,5){\line(-1,1){15}}
\put(0,5){\line(-1,0){15}}
\put(-15,20){\line(0,-1){15}}
\put(-2,10){\vector(1,0){10}}
\put(7,18){\framebox(4,4){15}}
\put(15,30){\line(1,0){5}}
\put(15,30){\line(0,1){5}}
\put(15,30){\line(-1,1){5}}
\put(15,30){\line(-1,0){5}}
\put(15,30){\line(-1,-1){5}}
\put(15,30){\line(0,-1){5}}
\put(15,30){\line(1,-1){5}}
\put(15,30){\line(1,1){5}}
\put(10,35){\line(0,-1){10}}
\put(10,35){\line(1,0){10}}
\put(20,25){\line(-1,0){10}}
\put(20,25){\line(0,1){10}}
\put(23, 28){\vector(2,-1){11}}
\put(7,-2){\framebox(4,4){14}}
\put(15,10){\line(1,0){5}}
\put(15,10){\line(0,1){5}}
\put(15,10){\line(-1,1){5}}
\put(15,10){\line(-1,0){5}}
\put(15,10){\line(-1,-1){5}}
\put(15,10){\line(0,-1){5}}
\put(15,10){\line(1,-1){5}}
\put(15,10){\line(2, -1){10}}
\put(25,5){\line(-1,1){10}}
\put(25,5){\line(-1,0){15}}
\put(10,15){\line(0,-1){10}}
\put(10,15){\line(1,0){5}}
\put(24,10){\vector(2,1){10}}
\put(23, 2){\vector(1,-1){9}}
\put(7, -31){\framebox(4,4){13}}
\put(15,-15){\line(1,0){5}}
\put(15,-15){\line(0,1){5}}
\put(15,-15){\line(-1,2){5}}
\put(15,-15){\line(-1,1){5}}
\put(15,-15){\line(-1,0){5}}
\put(15,-15){\line(0,-1){5}}
\put(15,-15){\line(-1,-2){5}}
\put(15,-15){\line(-1,-1){5}}
\put(10,-5){\line(0,-1){20}}
\put(20, -15){\line(-1,1){10}}
\put(20,-15){\line(-1,-1){10}}
\put(22,-15){\vector(1,0){10}}
\put(33,6){\framebox(4,4){12}}
\put(40, 17.5){\line(1,0){5}}
\put(40,17.5){\line(0,1){5}}
\put(45,17.5){\line(-1,1){5}}
\put(40,17.5){\line(-1,0){5}}
\put(40,17.5){\line(0,-1){5}}
\put(35,22.5){\line(1,-1){5}}
\put(40,17.5){\line(-1,-1){5}}
\put(40,17.5){\line(1,-1){5}}
\put(45,17.5){\line(-1,0){5}}
\put(35,12.5){\line(0,1){10}}
\put(45,17.5){\line(0,-1){5}}
\put(35, 12.5){\line(1,0){10}}
\put(40, 22.5){\line(-1,0){5}}
\put(43, 25){\vector(1,1){13}}
\put(47, 18){\vector(2,1){10}}
\put(47, 13){\vector(3, -1){12}}
\put(33, -26){\framebox(4,4){11}}
\put(40,-15){\line(1,0){5}}
\put(40,-15){\line(0,1){5}}
\put(40,-15){\line(-1,2){5}}
\put(40,-15){\line(-1,1){5}}
\put(40,-15){\line(-1,0){5}}
\put(40,-15){\line(0,-1){5}}
\put(40,-15){\line(-1,-1){5}}
\put(35,-5){\line(0,-1){15}}
\put(45, -15){\line(-1,1){10}}
\put(45,-15){\line(-1,-1){5}}
\put(35,-20){\line(1,0){5}}
\put(48,-15){\vector(1,0){10}}
\put(45,-10){\vector(1,1){12}}
\put(42, -8){\vector(1,2){14}}
\put(57, 33){\framebox(4,4){7}}
\put(65,40){\line(1,0){5}}
\put(65,40){\line(0,1){5}}
\put(65,40){\line(-1,1){5}}
\put(65,40){\line(-1,0){5}}
\put(65,40){\line(0,-1){5}}
\put(60,40){\line(1,-1){5}}
\put(65,40){\line(1,-1){5}}
\put(70,40){\line(-1,0){5}}
\put(60,40){\line(0,1){5}}
\put(65,45){\line(1,-1){5}}
\put(70,40){\line(0,-1){5}}
\put(65,35){\line(1,0){5}}
\put(65,45){\line(-1,0){5}}
\put(73, 35){\vector(1,-1){8}}
\put(57,14){\framebox(4,4){8}}
\put(65,25){\line(1,0){5}}
\put(65,25){\line(-1,0){5}}
\put(65,25){\line(-1,1){5}}
\put(60, 30){\line(2,-1){10}}
\put(65,25){\line(0,-1){5}}
\put(65,25){\line(1,-1){5}}
\put(65,25){\line(-1,-1){5}}
\put(60,20){\line(1,0){10}}
\put(70,25){\line(0,-1){5}}
\put(60,20){\line(0,1){10}}
\put(70,18){\vector(1,-1){10}}
\put(71,22){\vector(1,0){7}}
\put(57, -4){\framebox(4,4){9}}
\put(65,7){\line(1,0){5}}
\put(65,7){\line(0,1){5}}
\put(70,7){\line(-1,1){5}}
\put(65,7){\line(-1,0){5}}
\put(65,7){\line(0,-1){5}}
\put(65,7){\line(-1,-1){5}}
\put(70,7){\line(-1,-1){5}}
\put(65,7){\line(-1,1){5}}
\put(60,2){\line(0,1){10}}
\put(60,12){\line(1,0){5}}
\put(60, 2){\line(1,0){5}}
\put(69,3){\vector(2,-1){9}}
\put(70,10){\vector(1,1){8}}
\put(57, -26){\framebox(4,4){10}}
\put(65,-15){\line(1,0){5}}
\put(65,-15){\line(0,1){5}}
\put(65,-15){\line(-1,2){5}}
\put(65,-15){\line(-1,1){5}}
\put(65,-15){\line(-1,0){5}}
\put(65,-15){\line(-1,-1){5}}
\put(60, -5){\line(0,-1){15}}
\put(70, -15){\line(-1,1){10}}
\put(70,-15){\line(-2,-1){10}}
\put(70, -13){\vector(1,1){8}}
\put(79, 12){\framebox(4,4){5}}
\put(85,23){\line(1,0){5}}
\put(85,23){\line(0,1){5}}
\put(85,23){\line(-1,0){5}}
\put(85,23){\line(0,-1){5}}
\put(80,23){\line(0,-1){5}}
\put(80,18){\line(1,0){5}}
\put(90,23){\line(-1,-1){5}}
\put(85,28){\line(1,-1){5}}
\put(80,23){\line(1,1){5}}
\put(85,23){\line(-1,-1){5}}
\put(92,24){\vector(1,1){8}}
\put(90, 20){\vector(2,-1){9}}
\put(79, -11){\framebox(4,4){6}}
\put(85,0){\line(1,0){5}}
\put(85,0){\line(0,1){5}}
\put(85,0){\line(-1,1){5}}
\put(85,0){\line(-1,0){5}}
\put(85,0){\line(-1,-1){5}}
\put(80,5){\line(0,-1){10}}
\put(90, 0){\line(-1,1){5}}
\put(90,0){\line(-2,-1){10}}
\put(80,5){\line(1,0){5}}
\put(89, 4){\vector(2,1){9}}
\put(87, -5){\vector(2,-1){10}}
\put(99, 25){\framebox(4,4){2}}
\put(105,35){\line(1,0){5}}
\put(105,35){\line(0,1){5}}
\put(105,35){\line(-1,0){5}}
\put(105,35){\line(0,-1){5}}
\put(100,35){\line(1,-1){5}}
\put(105,40){\line(1,-1){5}}
\put(100,35){\line(1,1){5}}
\put(110,35){\line(-1,-1){5}}
\put(99,2){\framebox(4,4){3}}
\put(105,13){\line(1,0){5}}
\put(105,13){\line(0,1){5}}
\put(110,13){\line(-1,1){5}}
\put(105,13){\line(-1,0){5}}
\put(105,13){\line(-1,-1){5}}
\put(100,13){\line(0,-1){5}}
\put(110,13){\line(-2,-1){10}}
\put(105,18){\line(-1,-1){5}}
\put(113,13){\vector(1,0){10}}
\put(99, -21){\framebox(4,4){4}}
\put(105,-10){\line(1,0){5}}
\put(105,-10){\line(-1,0){5}}
\put(105,-10){\line(-1,-1){5}}
\put(105,-10){\line(-1,1){5}}
\put(100,-5){\line(0,-1){10}}
\put(110,-10){\line(-2,-1){10}}
\put(100,-5){\line(2,-1){10}}
\put(123,2){\framebox(4,4){1}}
\put(130,13){\line(1,0){5}}
\put(130,13){\line(0,1){5}}
\put(135,13){\line(-1,1){5}}
\put(130,13){\line(-1,-1){5}}
\put(135,13){\line(-2,-1){10}}
\put(130,18){\line(-1,-2){5}}
\end{picture}
\end{center}
\vspace{4.5cm}
\caption{Classification of nef toric surfaces.
The arrows indicate 
blow-downs. 
The numbers in frames are reference numbers.
Note that $S_3^0$, 
$S_4^0$, and  $S_5^0$ introduced in \S \ref{sec:deformation}
are No. 16, 14, and 12, respectively. }
\label{fig:blow-down}
\end{figure}
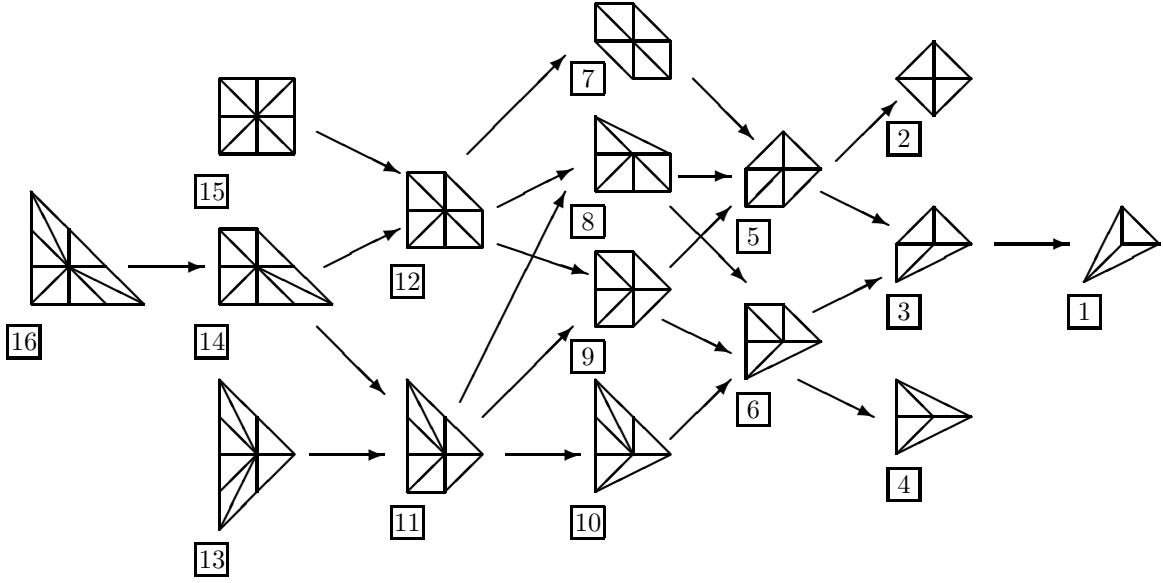

The following classification is due to Batyrev \cite{Batyrev}
(see also \cite[Table 1]{CKYZ}).
\begin{lemma}
There are exactly sixteen nef toric surfaces, whose fans are shown in Figure \ref{fig:blow-down}.  
\end{lemma} 

We will refer the nef toric surfaces using  the numbers
shown in frames 
in Figure \ref{fig:blow-down}.  

\begin{proof}
The minimal nef toric surfaces are $\Bbb{P}^2$, $\Bbb{P}^1 \times \Bbb{P}^1$, 
and the Hirzebruch surface $\Bbb{F}_2$,  
which are No. 1, No. 2, and No. 4 respectively. 
Nef toric surfaces are obtained from them 
by blowing up at a torus-fixed point successively. 
By the nef condition,  we must blow-up 
at a torus-fixed point which is not on a torus-fixed $(-2)$-curve. 
All possible patterns of blowing-ups are listed in Figure \ref{fig:blow-down}.  
Note that No. 13, 15, and 16 can no longer be blown-up to nef toric surfaces, 
since all of their torus-fixed points are on a 
torus-fixed $(-2)$-curve. This completes the classification. 
\end{proof}

\vspace{1cm}

\begin{proposition}\label{prop:def_type} 
A nef toric surface has a smooth versal deformation family of 
dimension $h^1(\Theta)$, whose general member is a del Pezzo surface of degree $c_1^2$.
\end{proposition}
$h^1(\Theta)$ and $c_1^2$ are given in Table \ref{tab:topo}. 

\begin{proof}

Note that $h^2(\Theta)=0$ for any smooth compact toric surface 
(Corollary \ref{cor:unob}). 
This implies smoothness of a versal deformation family \cite{KNS}. 

Versal deformation families of nef toric surfaces 
are constructed inductively as follows.
Let $\pi : \tilde{S} \to S$ be one of the blowing-ups in Figure \ref{fig:blow-down}.  
Let $P\in S$ be the center of the blowing-up $\pi$ which is the intersection of two torus-fixed curves 
$C_1$ and $C_2$ (see Figure \ref{fig:center}). 
By comparing Table \ref{tab:topo} with Figure \ref{fig:blow-down},  
we have
\begin{equation}\label{eq:moduli}
h^1(\tilde{S}, \Theta) = 
\begin{cases}
h^1(S, \Theta)  & \qquad \mathrm{if}~~ C_1^2 > -1,~ C_2^2 > -1~,\\  
h^1(S, \Theta) +1 & \qquad \mathrm{if}~~C_1^2 = -1,~ C_2^2 > -1~,\\
h^1(S, \Theta) +2 & \qquad \mathrm{if}~~C_1^2 = C_2^2 = -1~.
\end{cases}
\end{equation}
Since smooth rational curves on complex surfaces with self-intersection $\geq -1$ 
is stable under small deformations \cite[Example in p.86]{Kodaira} (see also \cite[IV. 3.1]{BPV}), 
a complete deformation family of $\tilde{S}$ can be found as a simultaneous blowing-up 
of a complete deformation family of $S$.  
Furthermore, by eq. (\ref{eq:moduli}), we can find a versal deformation family of $\tilde{S}$ as follows. 
First, we consider a versal deformation family $\mathcal{S}$
of $S$ on which $C_1$ and $C_2$ deform holomorphically. 
If both of $C_1$ and $C_2$ have self-intersection $> -1$, 
simultaneous blowing up of $\mathcal{S}$ at $P$ gives a versal deformation family of $\tilde{S}$
which is of dimension $h^1(S, \Theta)$.
If $C_1^2=-1$ and $C_2^2 >-1$, we move the center $P$ in the $C_2$ direction (see Figure \ref{fig:center}) and blow $\mathcal{S}$ up simultaneously  to get a versal 
deformation family of $\tilde{S}$ which is of dimension $h^1(S, \Theta) +1$. 
If $C_1^2 = C_2^2=-1$, we move the center $P$ 
in the whole direction 
and blow $\mathcal{S}$ up simultaneously  to get a versal deformation family of $\tilde{S}$
which is of dimension $h^1(S, \Theta)+2$. \\  
Thus we can find versal deformation families of nef toric surfaces inductively. 
It is easy to see that their general members are del Pezzo surfaces. 
\end{proof}
\begin{figure}[t]
\unitlength 0.1cm
\begin{center}
\begin{picture}(30,20)(20,-10)
\put(33, 13){\line(-1,-1){25}}
\put(10, 4){\line(1,0){40}}
\put(5, 4){$C_1$}
\put(35, 10){$C_2$} 
\put(24, 4){\circle*{2}}
\put(24, 0){$P$}
\thicklines
\put(24,4){\vector(1,1){5}}
\put(24,4){\vector(-1,-1){5}}
\end{picture}
\hspace*{3cm}
\begin{picture}(20,20)(20,-10)
\put(33, 13){\line(-1,-1){25}}
\put(10, 4){\line(1,0){40}}
\put(5, 4){$C_1$}
\put(35, 10){$C_2$} 
\put(24, 4){\circle*{2}}
\put(24, 0){$P$}
\thicklines
\put(23,6){\vector(-1,4){2}}
\put(22,6){\vector(-2,1){6}}
\end{picture}
\end{center}
\caption{The center $P$ of a blowing-up 
($C_1$ and $C_2$ are torus-fixed curves) and its moving.  The left is the case with 
$C_1^2=-1$, $C_2^2\geq 0$ and the right is the case with $C_1^2=C_2^2=-1$.}\label{fig:center}
\end{figure}
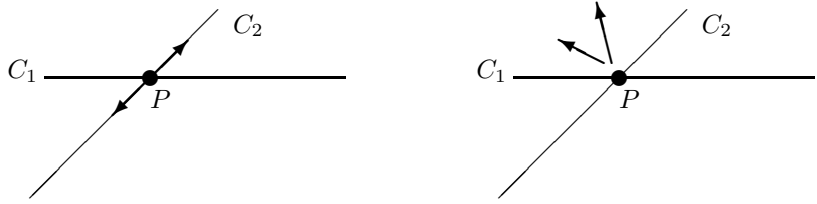
\begin{table}[ht]
\begin{equation}\notag
\begin{array}{|c||r|r|cr|cr|crrr|cr|crr|r|} \hline
\mathrm{Deformation~type}&I&II&III&&IV&&V&&&&VI&&VII&&&VIII \\ \hline
\mathrm{No.}            &1&3&2&4&5&6&7&8&9&10&11&12&13&14&15&16\\ \hline
c_1^2       &9&8&8&&7&&6&&&&5&&4&&&3               \\ \hline
c_2           &3&4&4&&5&&6&&&&7&&8&&&9      \\ \hline
-\frac{7c_1^2 - 5c_2}{6} &-8&-6&-6&&-4&&-2&&&&0&&2&&&4\\  \hline
h^0(\Theta)  &8&6&6&7&4&5&2&4&3&5&3&2&3&2&2&2       \\ \hline 
h^2(\Theta)  &0&0&0&0&0&0&0&0&0&0&0&0&0&0&0&0 \\ \hline
h^1(\Theta)  &0&0&0&1&0&1&0&2&1&3&3&2&5&4&4&6 \\ \hline
\end{array}
\end{equation}
\caption{Eight deformation types and 
$h^1(\Theta) \left(= - \frac{7c_1^2 - 5c_2}{6} + h^0(\Theta) + h^2(\Theta)\right)$~.}
\label{tab:topo}     
\end{table}
\section{Unobstructedness}\label{sec:unobstructed}
Let $X$ be a smooth compact toric surface,  
$D:=D_1+\cdots+D_r$  be the sum of all torus invariant divisors  $D_1, \cdots ,D_r$, 
and $\Theta(-\log D)$ be the sheaf of germs of holomorphic vector fields 
with logarithmic zeros along $D$.
\begin{lemma}\label{lem:toric}
$H^2(X,\Theta(-\log D)) = 0$. 
\end{lemma}
\begin{proof}
Since $\Theta(-\log D)=\mathcal{O} \otimes_{\mathbb{Z}} N$ 
(cf. \cite[Proposition 3.1]{Oda}),  
where $N$ is the $2$-dimensional lattice such that 
the fan of $X$ sits in $N\otimes \mathbb{R}$. 
$H^2(X,\Theta(-\log D))
 = H^2(X,\mathcal{O} \otimes_{\mathbb{Z}}N)
 = H^2(X,\mathcal{O} \oplus \mathcal{O}) = 0$, 
since $H^2(X,\mathcal{O})=0$ 
(cf. \cite[Corollary 2.8]{Oda}). 
\end{proof}

\begin{corollary}\label{cor:unob}
$H^2(X, \Theta) =0$. 
\end{corollary}
\begin{proof}
From the exact sequence (cf. \cite[Theorem 3.12]{Oda})
\begin{equation*}
\begin{CD}
0@>>>\Theta(-\log D)@>>>\Theta @>>>\oplus_{i=1}^{r} \mathcal{O}(D_i)|_{D_i}@>>>0~,
\end{CD}
\end{equation*}
and Lemma \ref{lem:toric}, we have $H^2(X, \Theta) = 0$.
\end{proof}

\section{$F_d$ $(1\leq d\leq 6)$}\label{sec:F_d}
Let $\mathsf{b[k]}:=(2\sin \frac{k\lambda}{2})^2$.
\begin{equation}\notag
F_1= \frac{1}{\mathsf{b[1]}} \,m(e_6)~,
\qquad
F_2= \frac{1}{2 \cdot \mathsf{b[2]}}\,m(2e_6)
+\frac{-2}{\mathsf{b[1]}}\,m(-e_1+l)~, 
\end{equation}
\begin{equation}\notag
\begin{split}
F_3&= \frac{1}{3\cdot \mathsf{b[3]}}\,m(3e_6) 
+\frac{3}{\mathsf{b[1]}}\,m(l) +
\Big(-4+\frac{27}{\mathsf{b[1]}}\Big)
\, m(-e_1-e_2-e_3-e_4-e_5-e_6+3l ) ~,
\end{split}
\end{equation}
\begin{equation}\notag
\begin{split}
F_4&=\frac{1}{4\cdot \mathsf{b[4]}}\,m(4e_6 ) 
+\frac{-2}{2\cdot \mathsf{b[2]}}\,m(-2 e_1+2l) 
+\frac{-4}{\mathsf{b[1]}}\,m(-e_1-e_2+2l) 
\\&+\Big(5+\frac{-32}{\mathsf{b[1]}}\Big)
\,m(-e_1-e_2-e_3-e_4-e_5+3l ) 
\end{split}
\end{equation}
\begin{equation}\notag
\begin{split}
F_5&=
\frac{1}{5\cdot \mathsf{b[5]}}\,m(5e_6) 
+\frac{5}{\mathsf{b[1]}}\,m(-e_1+2l)  
 +\Big(-6+\frac{35}{\mathsf{b[1]}}\Big)
  \,m(-e_1-e_2-e_3-e_4+3l) 
  \\
&+\Big(7\cdot \mathsf{b[1]}-68+\frac{205}{\mathsf{b[1]}}\Big)
   \,m(-2e_1-e_2-e_3-e_4-e_5-e_6+4l)~, 
\end{split}
\end{equation}
\begin{equation}\notag
\begin{split}
F_6&=\frac{1}{6\cdot \mathsf{b[6]}}\,m(6e_6)  
-\frac{2}{3\cdot \mathsf{b[3]}}\,m(-3e_1+3l) 
+\Big(\frac{3}{2\cdot \mathsf{b[2]}}-\frac{6}{\mathsf{b[1]}}\Big)\,m(2l) 
  \\
&+\Big(7-\frac{36}{\mathsf{b[1]}}\Big)
  \,m(-e_1-e_2-e_3+3l) 
  \\
&+ \Big(-8\cdot \mathsf{b[1]}+72 -\frac{198}{\mathsf{b[1]}}\Big)
   \,m(-2e_1-e_2-e_3-e_4-e_5+4l) 
  \\
&+\Big(9\cdot \mathsf{b[1]}^2-108\cdot \mathsf{b[1]}+498 -\frac{936}{\mathsf{b[1]}}\Big)
  \,m(-e_1-e_2-e_3-e_4-e_5-e_6+4l)
  \\
&+
 \Bigg(\frac{1}{2}\Big(-4+\frac{27}{\mathsf{b[2]}}\Big) 
  -10\cdot \mathsf{b[1]}^3+141\cdot \mathsf{b[2]}^2 -846\cdot \mathsf{b[1]}+2636  
-\frac{3780}{\mathsf{b[1]}}\Bigg)
\\&\times  m(-2e_1-2e_2-2e_3-2e_4-2e_5-2e_6+6l)~.
\end{split}
\end{equation}


\end{document}